\documentclass[onecolumn,final,12pt]{IEEEtran}

 \usepackage{graphicx}
\usepackage{amssymb}
\usepackage{amsmath}

\def\n{\noindent}
\def\cB{{\mathcal B}}
\def\cD{{\mathcal D}}
\def\cE{{\mathcal E}}

\def\cN{{\mathcal N}}

\def\RR{{\mathbb R}}

\def\EE{{\mathbb E}}
\def\Ex#1{{\EE\left\{#1\right\}}}

\def\tsigma{{\tilde\sigma}}

\def\tp{{\tilde p}}
\def\tepsilon{{\tilde \epsilon}}

\def\foorp{{\hfill$\spadesuit$}}

\newtheorem{lemma}{Lemma}

\newtheorem{definition}{Definition}
\newtheorem{theorem}{Theorem}

\begin{document}

%\begin{frontmatter}



\title{Determining local transientness of audio signals}


%
% elsevier
%\author[LATP]{St\'ephane Molla}
%\ead{molla@cmi.univ-mrs.fr}
%\author[LATP]{Bruno Torr\'esani}
%\ead{torresan@cmi.univ-mrs.fr}

%\address[LATP]{LATP, CMI, 39 rue Joliot-Curie,
%13453 Marseille cedex 13, France}

\author{St\'ephane Molla\thanks{Both authors are with
LATP, CMI, 39 rue Joliot-Curie, 13453 Marseille cedex 13, France;
email: molla@cmi.univ-mrs.fr; torresan@cmi.univ-mrs.fr.
Work supported in part by the European Union's
Human Potential Programme, under contract HPRN-CT-2002-00285 (HASSIP)}
and
Bruno Torr\'esani}

\maketitle

\begin{abstract}
We describe a new method for estimating the degree of
``transientness'' and ``tonality'' of a class of compound
signals involving simultaneously transient and harmonic features.
The key assumption is that both transient and tonal layers admit
sparse expansions, respectively in wavelet and local cosine bases.
The estimation
is performed using particular form of entropy (or theoretical dimension)
functions. We provide theoretical estimates on the behavior of the
proposed estimators, as well as numerical simulations.
Audio signal coding provides a natural field of application.
\end{abstract}

%\end{frontmatter}

%\newpage
% main text
\section{Introduction}
\label{se:intro}
Many generic signal classes feature significantly different
``components'', such as transients, (locally) sinusoidal or harmonic
``partials'', or stochastic-like components in sounds, or edges,
textures, etc. in images. Detecting the presence of such
components is one of the classical signal processing
problems (see for example ~\cite{Jaillet01detection}
and~\cite{Verma00perceptually} and references therein for reviews.) Another
interesting problem is to estimate whether a given portion of a signal
is for example more transient than harmonic or periodic,
or in other words to estimate ``transientness'' or ``tonality''
indices: quantitative measures of the local proportion of
transient and tonal features in a signal. Such indices find
immediate applications in several contexts, including the hybrid signal
coders~\cite{Daudet02hybrid,Levine99switched,Verma00perceptually}
which use different methods for encoding transient or tonal regions
(and were the main motivation of this work),
more general purpose hybrid models~\cite{Berger94removing}, or
similar recent ideas in image
coding~\cite{Meyer02multilayered,Romberg03Approximation}.
While there exist fairly standard tools for transient detection
or harmonic signal detection, the problem of quantitative measure
of proportion does not seem to have received much attention.

We propose here simple criteria, based on transform coding ideas,
for estimating such indices. The main idea is to use orthonormal bases
in signal spaces which are significantly different from each other
in the following sense: a given component has a sparse expansion
in a given basis, while the others have dense expansions. Information
theoretic criteria (we elaborate on the case of a variant of Shannon's
entropy) therefore yield estimates for the indices.

We focus here on the case of transient and locally sinusoidal
(or harmonic) layers in audio signals, using wavelet and local
cosine bases. However, the approach we develop may be adapted
to different signal layers (chirps for example), and to higher dimensions.
We provide theoretical estimates for the behavior of transientness
and tonality indices, and illustrate our results
by numerical simulations and tests on real sounds.

\section{Theoretical analysis}
\label{se:sparse}
We focus on the particular application to audio signals, and
limit ourselves to transient and tonal features. Our starting point
is the assumption that transient signals admit a sparse expansion
in a wavelet basis (provided the wavelets have small enough time support),
and that tonals admit a sparse expansion in local cosine basis (with
smooth enough window function.) We are naturally led to consider a
generic redundant ``dictionary'' made out of two such orthonormal
bases, denoted by $\psi_\lambda$ and $w_\delta$ respectively
(we refer to~\cite{Carmona98practical,Mallat98wavelet}
for detailed tutorials), and signal expansions of the form
\begin{equation}
\label{fo:sig.model}
x = \sum_{\lambda\in\Lambda} \alpha_{\lambda}\psi_{\lambda} +
\sum_{\delta\in\Delta} \beta_{\delta} w_{\delta} +r
= x_{tr} + x_{ton} +r\ ,
\end{equation}
where $\Lambda$ and $\Delta$ are (small, and this will be the main
{\em sparsity} assumption) subsets of the index sets,
termed {\em significance maps}.
The nonzero coefficients $\alpha_\lambda$
are assumed to be independent $\cN(0,\sigma_\lambda)$
random variables, and the nonzero coefficients $\beta_\delta$
are assumed to be independent $\cN(0,\tsigma_\delta)$ random variables.
$r$ is a residual signal, which is not sparse with respect
to the two considered bases (we shall talk of {\em spread residual}),
and is to be neglected or described differently.

Given a signal assumed for simplicity to be of the
form~(\ref{fo:sig.model}), with unknown values of $|\Delta|$
and $|\Lambda|$ (the cardinalities of $\Delta$ and $\Lambda$
respectively), we are interested in finding estimates for the latter.
More precisely, we seek estimates for the ``transientness'' and
``tonality'' indices
\begin{equation}
\label{fo:indices}
I_{tr} = \frac{|\Lambda|}{|\Delta|+|\Lambda|}\ ,\qquad
I_{ton} = \frac{|\Delta|}{|\Delta|+|\Lambda|}\ .
\end{equation}
We propose a procedure close to the notions of {\em theoretical
dimension} or $\alpha$-entropies, advocated by
Wickerhauser (see~\cite{Wickerhauser94adapted} for a review.)
%, which may in some situations be shown to be closely
%connected to the notion of Shannon entropy~\cite{Trgo95relation}.
Our approach is based upon the following heuristics. Consider a
signal $x$, and expand it into an orthonormal basis. Estimating
the ``size'' of $x$ in this basis may clearly be done by counting
the number of nonzero coefficients (the $\ell^0$ norm of the
sequence of coefficients), or the number of coefficients above some
threshold. It has been shown~\cite{Wickerhauser94adapted}
that alternative approaches are possible, including $\ell^p$ norms
(with $p<2$) or entropy, yielding comparable results.
Considering a hybrid signal as in~(\ref{fo:sig.model}) neither
its wavelet expansion nor its local cosine expansion will be sparse.
However, since by assumption only very few coefficients
$\alpha_\lambda$ and $\beta_\delta$
are nonzero, most wavelet coefficients
$\langle x,\psi_\lambda\rangle$ actually originate from
the tonal part $x_{ton}$, and most local cosine coefficients
$\langle x,w_\delta\rangle$ originate from the transient part $x_{tr}$.
Therefore, calculating these $\ell^p$ norms or entropy
%one of the above mentioned  quantities
from the wavelet coefficients $\langle x,\psi_\lambda\rangle$
is expected to provide (approximately) estimates on the number of
nonzero (or significant) $\beta_\delta$ coefficients, and vice versa.
We elaborate below on the specific case of the logarithmic
dimensions, for which such a behavior may be proved.
For the sake of simplicity, we shall work in this section
in a finite dimensional context.
\begin{definition}
Given an orthonormal basis
$\cB = \{e_n,n=1,\dots N\}$
of a given $N$-dimensional signal space $\cE\cong \RR^N$,
define the logarithmic dimension of
$x\in\cE$ in the basis $\cB$ by
\begin{equation}
\cD_\cB(x) = \frac1{N}\,\sum_{n=1}^N
\log_2\left(|\langle x,e_n\rangle|^2\ \right)\ .
\end{equation}
\end{definition}
%In a nutshell, for a fixed $x\in\cE$, the smaller $\cD_\cB(x)$, the
%sparser the representation of $x$ in the basis $\cB$. In addition,
%$2^{\cD_\cB(x)}$ provides a rough estimate of the number of
%significant coefficients to represent $x$ in the basis $\cB$.

It follows from a simple calculation that in the framework of
the signal models under consideration, 
\begin{lemma}
Given an orthonormal basis $\cB = \{e_n,n=1,\dots N\}$,
assuming that the coefficients $\langle x,e_n\rangle$
are $\cN(0,\sigma_n)$ random variables, one has
\begin{equation}
\Ex{\cD_\cB(x)} = %\left(1+\frac\gamma{\log(2)}\right) 
C+ \frac1{N}\sum_{n=1}^N
\log_2(\sigma_n^2)\ ,
\end{equation}
where $C= 1+\gamma/\ln(2)\approx 1.832746177$ is a universal constant
($\gamma\approx .5772156649$ being Euler's constant.)
\end{lemma}

Returning to the model~(\ref{fo:sig.model}), and assuming that
the coefficients $\alpha_\lambda,\lambda\in\Lambda$
and $\beta_\delta,\delta\in\Delta$ are respectively $\cN(0,\sigma_\lambda)$
and $\cN(0,\tsigma_\delta)$ independent random variables, the coefficients
$a_\lambda = \langle x,\psi_\lambda\rangle$, 
$b_\delta = \langle x,w_\delta\rangle$,
are zero-mean normal random variables, whose variance depends
on whether $\lambda\in\Lambda$ (or $\delta\in\Delta$) or not.
For example,
\begin{equation}
\hbox{var}\{a_\lambda\} = \left\{
\begin{array}{ll}
\sigma_\lambda^2 + \sum_{\delta \in\Delta}\tsigma_\delta^2
|\langle \psi_\lambda,w_\delta\rangle|^2 &\hbox{ if }\lambda\in\Lambda\\
\sum_{\delta \in\Delta}\tsigma_\delta^2
|\langle \psi_\lambda,w_\delta\rangle|^2 &\hbox{ if }\lambda\not\in\Lambda\\
\end{array}
\right.
\end{equation}
and we obtain, for the $\Psi=\{\psi_\lambda\}$ basis
\begin{equation}
\label{fo:log.dim.1}
\begin{split}
\Ex{\cD_\Psi(x)}\! =\! C\! &\!+\! 
\frac1{N}\!\log_2\!\!
\bigg[\!\prod_{\lambda\in\Lambda}\!\!\left(\!\sigma_\lambda^2 +\!\!
\sum_{\delta\in\Delta}
\tsigma_\delta^2 |\langle \psi_\lambda,w_\delta\rangle|^2\!\right)\\
&\times
\prod_{\lambda'\not\in\Lambda}\!\!\left(\sum_{\delta\in\Delta} \tsigma_\delta^2
|\langle \psi_{\lambda'},w_\delta\rangle|^2\right)\bigg]\ ,
\end{split}
\end{equation}
and a similar expression for the logarithmic dimension
$\cD_W(x)$ with respect to the $W=\{w_\delta\}$ basis.

In the simpler case where $\sigma_\lambda =\sigma$,
$\forall\lambda\in\Lambda$ and $\tsigma_\delta =\tsigma$,
$\forall\delta\in\Delta$, we introduce the Parseval weights
\begin{equation}
p_\lambda(\Delta) = \sum_{\delta\in\Delta}
|\langle w_\delta,\psi_{\lambda}\rangle|^2\ ,
\quad
\tp_\delta(\Lambda) = \sum_{\lambda\in\Lambda}
|\langle w_\delta,\psi_{\lambda}\rangle|^2\ .
\end{equation}
The following property is an immediate
consequence of Parseval's formula
(i.e. for all $f$, 
$\sum_\lambda |\langle f,\psi_\lambda\rangle|^2 = \|f\|^2$.)
\begin{lemma}
\label{le:parseval}
The Parseval weights satisfy
$$
0\le p_\lambda(\Delta)\le 1\ ,\quad
0\le \tp_\delta(\Lambda)\le 1\ .
$$
\end{lemma}
Introducing the {\em relative redundancies} of the bases
$\Psi$ and $W$ with respect to the significance maps
\begin{equation}
\epsilon(\Delta) = \sup_{\lambda\in\Lambda} p_\lambda(\Delta)\ ,
\quad
\tepsilon(\Lambda) = \sup_{\delta\in\Delta} p_\delta(\Lambda)\ ,
\end{equation}
we obtain simple estimates for the logarithmic dimension.
\begin{theorem}
With the above notations, assuming that the %significant
coefficients $\{\alpha_\lambda,\lambda\in\Lambda\}$ and
$\{\beta_\delta,\delta\in\Delta\}$ are independent identically distributed
$\cN(0,\sigma)$ and $\cN(0,\tsigma)$ normal variables respectively,
and assuming $r=0$, the following bounds hold
\begin{equation}
\label{fo:sup.bound}
\begin{split}
\quad\Ex{\cD_\Psi(x)} \ge C &
+ \frac{|\Lambda|}N\log_2(\sigma^2)\\
& 
+
\log_2\left(\prod_{\lambda'\not\in\Lambda}\! \left(
\tsigma^2 p_{\lambda'}(\Delta)\right)^{1/N}\right)
\end{split}
\end{equation}
\begin{equation}
\label{fo:inf.bound}
\begin{split}
\quad\Ex{\cD_\Psi(x)} \le C &+ \frac{|\Lambda|}N
\log_2(\sigma^2 +\epsilon(\Delta)\tsigma^2)\\
& +
\log_2\left(\prod_{\lambda'\not\in\Lambda}\! \left(
\tsigma^2 p_{\lambda'}(\Delta)\right)^{1/N}\right)\ ,
\end{split}
\end{equation}
with $C= 1+\gamma/\ln(2)\approx 1.832746177$.
Exchanging the roles of $\Delta$ and $\Lambda$, a similar bound
holds for the other logarithmic dimension $\cD_W(x)$.
\end{theorem}
\n{\em Proof: } The proposition follows directly from the fact that
in such a situation, equation~(\ref{fo:log.dim.1}) reduces to
\begin{equation}
\label{fo:log.dim.2}
\begin{split}
\Ex{\cD_\Psi(x)} = C + 
&\log_2\bigg(\prod_{\lambda\in\Lambda}\! \left(\sigma^2 +
\tsigma^2 p_\lambda(\Delta)\right)^{1/N}\\
&\times\prod_{\lambda'\not\in\Lambda}\! \left(
\tsigma^2 p_{\lambda'}(\Delta)\right)^{1/N}\bigg)\ ,
\end{split}
\end{equation}
from Lemma~\ref{le:parseval} and the definition of
$\epsilon(\Delta)$.\foorp

\medskip
This result may be understood and utilized as follows.
First notice that the bounds in Equations~(\ref{fo:sup.bound})
and~(\ref{fo:inf.bound}) differ by
$|\Lambda|\log_2(1+\epsilon(\Delta)\tsigma^2/\sigma^2)/N$.
Let us assume for a while that this term may be neglected
(more on that below.) Then the behavior of $\Ex{\cD_\Psi(x)}$
is essentially controlled by
$\log_2\left(\prod_{\lambda'\not\in\Lambda}\! \left(
\tsigma^2 p_{\lambda'}(\Delta)\right)^{1/N}\right)$.
The behavior of this term is not easy to understand, but a first
idea may be obtained by replacing $p_{\lambda'}(\Delta)$ by its
``ensemble average''
$
\frac1{N}\sum_{\lambda=1}^N p_{\lambda}(\Delta) = 
\frac1{N}\sum_{\delta\in\Delta}\|w_\delta\|^2 =
\frac{|\Delta|}N\ ,
$
which yields the approximate expression:
\begin{equation}
\label{fo:app}
\begin{split}
\Ex{\cD_\Psi(x)} \approx C
&
+ \frac{|\Lambda|}N \log_2(\sigma^2)
\\
&
+ \left(1-\frac{|\Lambda|}N\right)\,
\log_2\left(\tsigma^2\frac{|\Delta|}N\right)\ .
\end{split}
\end{equation}
Therefore, if the ``$\Psi$-component'' of the signal is sparse enough,
i.e. if $|\Lambda|/N$ is sufficiently small (compared with 1),
$\Ex{\cD_\Psi(x)}$ may be expected to behave as 
$\log_2\left(\tsigma^2\frac{|\Delta|}N\right)$. Set 
\begin{equation}
\label {fo:Npsi}
\hat N_{\Psi}(x) = 2^{\cD_{\Psi}(x)}\ .
\end{equation}
Replacing $\cD_{\Psi}(x)$ with its expectation, we see that
$\hat N_{\Psi}(x)\approx 2^C\,\tsigma^2\frac{|\Delta|}N$, which
yields an estimate (up to the multiplicative constant $2^C\tsigma^2/N$)
for the ``size'' of the tonal component of the signal.
Similarly, defining
\begin{equation}
\label {fo:NW}
\hat N_{W}(x) = 2^{\cD_{W}(x)}
%\approx 2^C\,\sigma^2\frac{|\Lambda|}N
\end{equation}
we obtain a similar estimate (up to the multiplicative constant
$2^C\sigma^2/N$) for the ``size'' of the $\Psi$
component of the signal. Both $\hat N_{\Psi}(x)$ and
$\hat N_{W}(x)$ are computable from the signals wavelet and
local cosine expansions, and we finally consider their
relative proportions, or ``rates''
\begin{equation}
\label{fo:est.rates}
\hat I_{ton} = \frac{\hat N_\Psi}{\hat N_\Psi +\hat N_W}\ ,\quad
\hat I_{tr} = \frac{\hat N_W}{\hat N_\Psi +\hat  N_W}\ ,
\end{equation}
which provide the desired estimates for the indices in Eq.~(\ref{fo:indices}).

\medskip
A few comments are in order here.

\begin{itemize}
\item[{\em i.}] The difference between the lower and
upper bounds depends on the sparsity
$|\Lambda|/N$ of the $\Psi$-component and the relative
redundancy parameters $\epsilon(\Delta)$. The latter actually describe
the intrinsic differences between the two considered bases.
When the bases are significantly different, the relative redundancy
may be expected to be small (notice that in any case, it is smaller than 1.)
\item[{\em ii.}] The relative redundancy parameters $\epsilon$
and $\tepsilon$ %which pop up in our model
differ from the one which is generally considered in the literature,
namely the {\em coherence} 
$
M[W\cup\Psi]
%= \sup_{\stackrel{b,b'\in W\cup\Psi}{b\ne b'}}|\langle b,b'\rangle|\ .
= \sup_{b\in W,b'\in \Psi}|\langle b,b'\rangle|
$
of the dictionary $W\cup\Psi$ (see
e.g.~\cite{Donoho01uncertainty,Elad02generalized,Gribonval03sparse}.)
The latter is intrinsic to the dictionary, while the Parseval weights
and corresponding $\epsilon$ and $\tepsilon$
provide a finer information, as they also account for the
signal models, via their dependence in the significance maps $\Lambda$ and
$\Delta$.
\item[{\em iii.}] Precise estimates for $\epsilon$ and
$\tepsilon$ are difficult to obtain (numerical simulations
yield values around 1/4.) More precise models
for the significance maps $\Delta$ and $\Lambda$ could provide better
understanding. In particular, structured models such as those
described in~\cite{Molla03hybrid} (implementing time persistence
in $\Delta$ and scale persistence in $\Lambda$)
are expected to yield smaller values for the relative redundancies than models
featuring uniformly distributed significance maps.
\end{itemize}

%\bigskip
Another interesting point is the sensitivity of such tools with respect
to departures from the model, or noise. We show that results similar to
the above ones still hold true in the presence of white noise,
i.e. assuming that the residual $r$ in~(\ref{fo:sig.model}) is
a zero-mean Gaussian white noise.
In such a situation, denoting by $s^2$ the variance of the
noise $r$, equation~(\ref{fo:log.dim.1}) becomes
\begin{equation*}
\label{fo:log.dim.1bis}
\begin{split}
\Ex{\cD_\Psi(x)}\! =\! C &\!+ \!
\frac1{N}\log_2\!\!
\bigg[\prod_{\lambda\in\Lambda}\!\!\left(\!\sigma_\lambda^2\! +\!\!
\sum_{\delta\in\Delta}\!\!
\tsigma_\delta^2 |\langle \psi_\lambda,w_\delta\rangle|^2\! +\! s^2\!\right)\\
&\times
\prod_{\lambda'\not\in\Lambda}\!\!\left(\sum_{\delta\in\Delta} \tsigma_\delta^2
|\langle \psi_{\lambda'},w_\delta\rangle|^2+s^2\right)\bigg]\ ,
\end{split}
\end{equation*}
and a similar expression for the logarithmic dimension
$\cD_W(x)$ with respect to the $W=\{w_\delta\}$ basis.
Hence, the approximate expression~(\ref{fo:app}) becomes
\begin{equation*}
\label{fo:app.bis}
\begin{split}
\Ex{\cD_\Psi(x)} \approx C
&
+ \frac{|\Lambda|}N \log_2(\sigma^2 + s^2)
\\
&
+ \left(1-\frac{|\Lambda|}N\right)\,
\log_2\left(\tsigma^2\frac{|\Delta|}N + s^2\right)\ .
\end{split}
\end{equation*}
The discussion above (suitably adaptated) still holds 
as long as the signal energy $\tsigma^2|\Delta|$ exceeds the noise
energy $s^2N$.

These estimates yield the following algorithm for estimating
transientness and tonality indices for sparse hybrid signals
$x$ as given in~(\ref{fo:sig.model}).
\begin{enumerate}
\item Compute logarithmic dimensions $\cD_\Psi(x)$ and $\cD_W(x)$
\item Compute $\hat N_\Psi(x)$ and $\hat N_W(x)$ as
in~(\ref{fo:Npsi}) and~(\ref{fo:NW}).
\item Compute the estimated rates as in~(\ref{fo:est.rates}).
\end{enumerate}

\section{Numerical results}
\label{se:num}
\begin{figure}
\centerline{
\includegraphics[width=9cm,height=5cm]{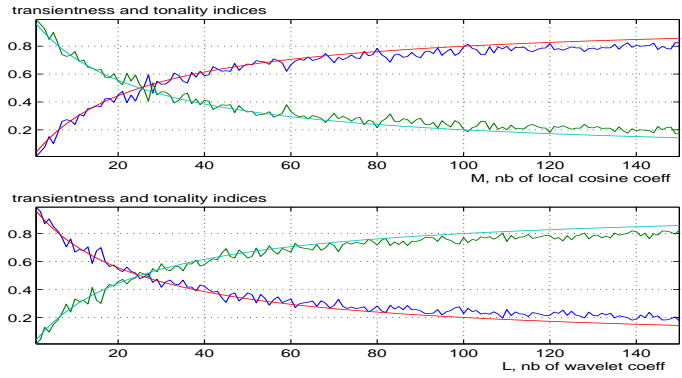}
}
\caption{Transientness and tonality estimates for the model (averaged
over 10 realizations.)
{\em Top:} $L=25$, and $M\in\{1,\dots 150\}$; increasing curves: $I_{ton}$
and $\hat I_{ton}$; decreasing curves: $I_{tr}$ and $\hat I_{tr}$.
{\em Bottom:} $M=25$, and $L\in\{1,\dots 150\}$; increasing curves: $I_{tr}$
and $\hat I_{tr}$; decreasing curves: $I_{ton}$ and $\hat I_{ton}$.
}
\label{fi:25}
\end{figure}

\begin{figure}
\centerline{
\includegraphics[width=9cm,height=5cm]{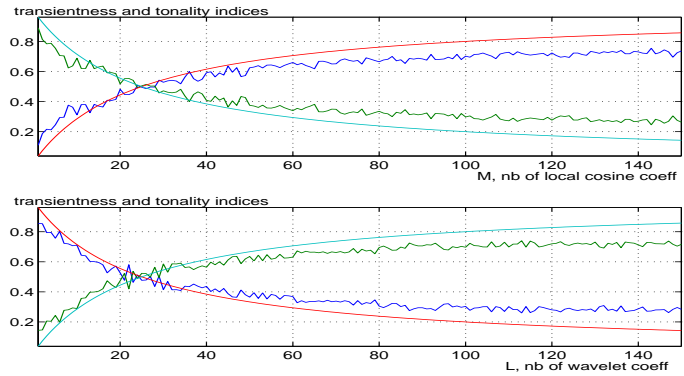}
}
\caption{Influence of white noise:
Transientness and tonality estimates for the model (averaged
over 10 realizations) with additional white noise. Same legends as before.
}
\label{fi:25n3}
\end{figure}

We generated several realizations of the signal model (with
$r=0$ first), with variable numbers $L$ of wavelet atoms and fixed number
$M$ of local cosines and vice versa,
and computed the estimated rates
$\hat I_{ton}$ and $\hat I_{tr}$,
to be compared with the ground truth~(\ref{fo:indices}), i.e.
$I_{ton} = M/(M+L)$ and $I_{tr} = L/(M+L) =1-I_{ton}$.
As may be seen from %the numerical simulations presented in
Figure~\ref{fi:25} (which corresponds to averages
over 10 realizations of the model), the estimated curves reproduce
quite well the correct ones. Some discrepancies may be observed
at the right hand side of the curves, where the sparsity assumptions
are not valid any more, and the correction terms in~(\ref{fo:app})
come into play. Observe that the curves cross precisely at the
correct location $M=L$.
The influence of the noise may be seen in Figure~\ref{fi:25n3}:
a white noise, whose energy equals 30\% of the signal's energy,
has been added. The effect is what can be anticipated
from~(\ref{fo:app.bis}), namely the presence of an additional noise
term moves the experimental curves away from the theoretical ones.

%\bigskip

\begin{figure}
\centerline{
\includegraphics[width=9cm,height=4.95cm]{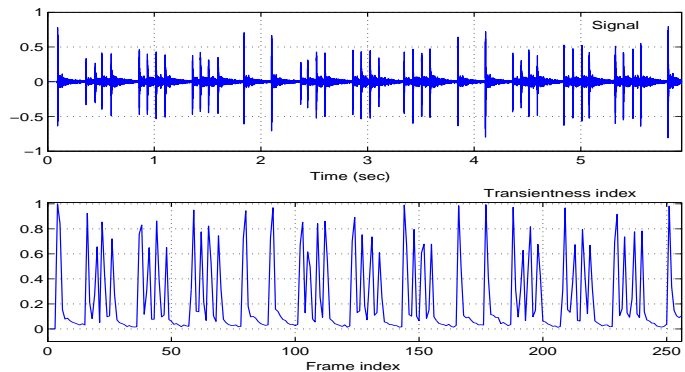}
}
\caption{Transientness index for the test 'castanet' signal.
Signal (top) and transientness index (bottom.)
}
\label{fi:casta}
\end{figure}

Besides the numerical simulations above, the transientness and
tonality indices have been tested on real audio signals, yielding
very sensible results.\footnote{Additional material, including
sound files, may be found at the web site
{\tt http://www.cmi.univ-mrs.fr/\~{ }torresan/papers/balance}.}
A first example, based upon a simple 'castanet' signal (6 s long,
sampled at 44,100 kHz) is shown in Figure~\ref{fi:casta}. A value for the
transientness index and the tonality index was computed for all time
frames (23 ms long.) Since $I_{ton}=1-I_{tr}$, only the transientness
index is displayed for the sake of clarity. This signal is quite simple, as it
essentially exhibits attacks followed by harmonic tones, and is thus
a ``perfect'' test for the proposed approach. As may be seen from
the bottom plot of Figure~\ref{fi:casta}, all attacks are
correctly captured, and the corresponding index is quite
high. In between attacks, the transientness index is very low,
which is also natural since the signal is essentially harmonic, thus
sparsely represented by local cosine basis.

\begin{figure}
\centerline{
\includegraphics[width=9cm,height=4.95cm]{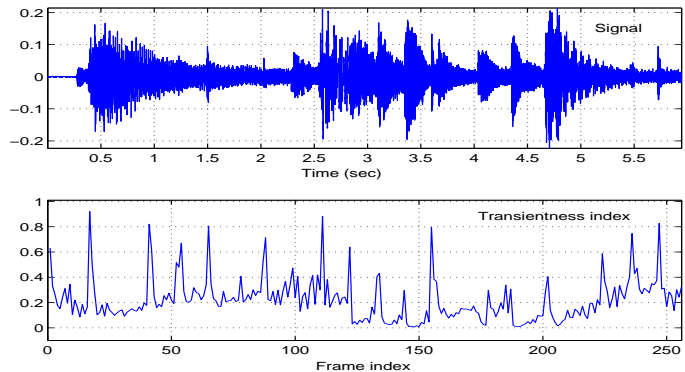}
}
\caption{Transientness index for the test ``jazz'' signal.
Signal (top) and transientness index (bottom.)
}
\label{fi:you}
\end{figure}

The second sound example displayed here is a more complex audio
signal, extracted from a jazz recording (about 6 s. long,
sampled at 44,100 kHz)
which features ``mixed'' tonals and transients.
The numerical results are displayed  in Figure~\ref{fi:you}.
Notice again that the ``obvious'' attacks of the signal have been
captured by the method. A closer examination of the signal (using a
``spectrogram type'' representation, not shown here) shows that
in the middle part of the signal (more precisely, between seconds
3 and 5), the harmonic content is stronger,
which explains the lower average value of $I_{tr}$ there.
This illustrates the fact that $I_{tr}$ really provides an estimate of
the {\em proportion} of transients relative to tonals, rather than an
absolute indicator of the presence of transient, such as the ones used
in transient detection~\cite{Levine99switched} for
example.\footnote{
\tt http://www.cmi.univ-mrs.fr/\~{ }torresan/papers/balance.}

More numerical results, in the framework of the hybrid audio coding
scheme developed in~\cite{Daudet02hybrid}, will be given in a
forthcoming publication~\cite{Molla03hybrid}.

\section{Conclusions}
\label{se:concl}
Sparsity of wavelet and local cosine signal representations
may be exploited in order to estimate the relative amount of tonal and
transient components present in the signal. This approach proves to be
extremely effective in the context of hybrid audio signal
coding~\cite{Daudet02hybrid,Molla03hidden}, and possesses a wider
range of applications, including image coding~\cite{Meyer02multilayered}.

The theoretical analysis we have presented is based on strong
a priori assumptions on the signal (essentially, a hybrid model such
as~(\ref{fo:sig.model}), with sparse significance maps $\Lambda$
and $\Delta$.)
While this sparsity assumption is completely necessary, the equality
of variances may be relaxed; in that situation, the indices provide
estimates on the proportion of energies
$\sigma^2 |\Lambda|$ and $\tsigma^2 |\Delta|$ of the
two layers, rather than their size $|\Lambda|$ and $|\Delta|$.
This is the case for the numerical results on real signals, for
which the variances of the two layers are not known.

% and equal (or comparable) energies for the two layers.
%When this is not the case, the approach may easily be refined to account
%for departures from such a situation.

Finally, let us simply mention that the approach may be extended to
more than two layers, provided that the considered orthonormal bases
are sufficiently different (in terms of their ``Parseval weights'',
see above) to allow the separation.  Again, this may prove useful
in the context of image coding, where new types of waveforms (e.g. curvelets)
may be introduced.

% \appendix

% \section{}
% \label{}

\end{document}